\documentstyle[12pt, amstex,epsf, epic, eepic]{amsart}
\renewcommand{\pf}{\noindent{\em Proof: }}

\newcommand{\s}[1]{\mathcal{#1}}

\newcommand{\Kernel}{\operatorname{Ker}}

\newcommand{\rank}{\operatorname{rank}}

\newcommand{\End}{\operatorname{End}}
\newcommand{\corank}{\operatorname{corank}}

\newcommand{\pid}{\operatorname{p.i. deg}}

\newtheorem{Thm}{Theorem}[section]
\newtheorem{Def}[Thm]{Definition} \newtheorem{Rem}[Thm]{Remark}
\newtheorem{Lem}[Thm]{Lemma} \newtheorem{Cor}[Thm]{Corollary}
\newtheorem{Prop}[Thm]{Proposition}

\numberwithin{equation}{section}

\begin{document}

\fontsize{.4cm}{.4cm}\selectfont\sf

\title[Quantized rank $R$ matrices]{Quantized rank $R$ matrices}
\author[Jakobsen and J\o ndrup]{Hans Plesner Jakobsen and S\o ren J\o ndrup} 
\address{
  Mathematics Institute\\ Universitetsparken 5\\ DK--2100 Copenhagen {\O},
  Denmark}\email{jakobsen@@math.ku.dk, jondrup@@math.ku.dk}

\date{\today}

\begin{abstract} First some old as well as new results about
  P.I. algebras, Ore extensions, and degrees are presented. Then
  quantized $n\times r$ matrices as well as certain quantized factor
  algebras $M^{r+1}_q(n)$ of $M_q(n)$ are analyzed. For
  $r=1,\dots,n-1$, $M^{r+1}_q(n)$ is the quantized function algebra of
  rank $r$ matrices obtained by working modulo the ideal generated by
  all $(r+1)\times (r+1)$ quantum subdeterminants and a certain
  localization of this algebra is proved to be isomorphic to a more
  manageable one. In all cases, the quantum parameter is a primitive
  $m$th roots of unity. The degrees and centers of the algebras are
  determined when $m$ is a prime and the general structure is obtained
  for arbitrary $m$.\footnotemark\footnotemark
\end{abstract}

\maketitle
\footnotetext[1]{MSC-class: 16S36, 17B37, 17B10.}  \footnotetext[2]{Key
  words: Quantized function algebra, quantum minor, factor algebra,
  P.I. algebra, block diagonal.}

\section{Introduction}

Through the last several years, quantized function algebras have
attracted a lot of attention (\cite{frt}, \cite{cp}, \cite{cl},
\cite{dd}, \cite{drin}, \cite{ho-la}, \cite{jaz1}, \cite{lev},
\cite{leso}, \cite{lu1}, \cite{manin}, \cite{pw}, and many others).
Amongst these, $M_q(n)$ has attracted the most attention. Since in
fact a number of candidates for the quantized function algebra of
$n\times n$ matrices have been proposed, we stress that the one we
consider here is the ``original'' (or ``standard'' or ``official'')
one introduced by Faddeev, Reshetikhin, and Takhtajan in \cite{frt}.

We wish to consider some natural subalgebras and quotients of this
algebra, namely the subalgebra $M_q(n,r)$ of quantized $n\times r$
matrices,  the subalgebra $A_{n,r}$ obtained by removing the 
$(n-r)\times(n-r)$ corner generated by those $Z_{i,j}$ for which
$i,j\geq r+1$, and finally the quotients
$M^{r+1}_q(n)=M_q(n)/I^{r+1}_q$ obtained by factoring out the ideal
$I^{r+1}_q$ generated by all $(r+1)\times (r+1)$ quantum
subdeterminants. The emphasis will throughout be on the case where $q$
is a primitive root of unity.

The major tool is the theory developed by De Concini and Procesi in
\cite{cp} as well as the theory of P.I. algebras. We have found it
convenient to collect these results, some corollaries to them, as well 
as some further developments in Section~\ref{pi} following immediately
after this introduction. In some sense, the results of  De Concini and Procesi  turn the
problem into an elementary one, but which at the same time is of a
such kind that one should not expect general results except possibly
in special cases.  Indeed, a major part of the procedure is to bring
into block diagonal form an integer coefficient skewsymmetric form. 

The new results we present relate to (iterated) skew polynomial
extensions and are particularly useful for the algebras $A_{n,r}$.
Even for the known case (c.f. \cite{jaz1}) $M_q(n)$ they are
sufficient, and we have chosen to illustrate this in
Section~\ref{stand}. Actually, the case of $M_q(n)$ was brought to a
completion by the discovery of a very special phenomenon for the
associated quasipolynomial algebra (\cite{jaz1}) and a substantial
further development of this observation now makes it possible to
attack $M_q(n,r)$. This is done in Section~\ref{nr}.

Section~\ref{qm} is devoted to the proof of the isomorphism
$A_{n,r}[d^{-1}]\simeq M^{r+1}_q(n)[d^{-1}]$, where $d$ is the quantum
$r\times r$ minor corresponding to the subalgebra generated by those
$Z_{i,j}$ for which $1\leq i,j\leq r$. A major tool is the
representation theory of quantized enveloping algebras. Having
established that, we turn our attention to $A_{n,r}$ in
Section~\ref{qfa}. Indeed, for questions relating to degree, center,
etc., it is sufficient to consider this algebra, which is more
manageable. The methods that worked well for $M_q(n,r)$ do not apply
as easily but, fortunately, the results obtained in Section~\ref{pi}
are applicable, especially after some fortunate guesses relating to
the center. As one consequence we obtain (combine
Corollary~\ref{com-cor} with Theeorem~\ref{com-thm}): {\em If $q$ is an
  odd primitive $m$th root of unity then $\deg M^{r+1}_q(n)
  =m^{nr-r(r+1)/2}$}.
 
We thank Ken  Goodearl and Tom Lenagan for criticizing the original
proof of Proposition~\ref{in-pro}.

\medskip

\section{P.I. algebras and their degrees}
\label{pi}
In this section, unless explicitly stated, $A$ denotes a prime P.I.
algebra and $k$ an algebraically closed field $k$ of characteristic
$0$. We assume throughout that $A$ is finitely generated (affine) as
an algebra over $k$.

We start by recalling some basic results from the theory of P.I.
algebras and then we show these results can be applied to calculations
of the degree of an algebra. Let us first recall some basic definitions.

\medskip

\begin{Def}[De Concini-Procesi, {\cite[p. 50] {cp}}] An algebra $A$ is said to be of degree at most $d$, if $A$ satisfies all identities of $(d\times d)$-matrices over a commutative  ring. If no such $d$ exists $A$ is said to have infinite degree. In either case the smallest possible $d$ is denoted by $\deg A$.
\end{Def}

\begin{Def}The P.I. degree of an algebra $A$, $\pid A$, is $\frac b2$, where $b$ is the smallest possible degree
of a multilinear polynomial which vanishes on $A$.
\end{Def}

If $A$ is a prime P.I.-algebra, such as ours,  one gets

\begin{Prop}[McConnel-Robson, {\cite[13.6.7 (v)]{mr}}] 
  $\deg A$ =  $\pid A$.
\end{Prop}

For an affine prime P.I. algebra $A$ one has several useful results
concerning the (P.I.)-degree. To state and prove these, let us first
recall that the intersection of all primitive ideals of $A$ is $0$
\cite[Theorem 13.10.3]{mr} and all primitive ideals are
maximal by Kaplansky's Theorem \cite[Theorem 13.3.8]{mr}.

\medskip

\begin{Prop}
 $\pid A=\sup_M\pid A/M$,  where $M$ runs through the
set of all maximal two-sided ideals.
\end{Prop}

\medskip

\pf Since $A/M$ is a factor algebra of $A$ for all
maximal ideals $M$, we get that the right hand side is bounded by the
left hand side.

From \cite[Corollary 13.6.7]{mr} we get that $\pid A/M\leq n$  if and only if $S_{2n}$, the standard identity, is an
identity for $A/M$, thus $\sup_M\pid A/M\leq n$ implies that
$S_{2n}$ is an identity for $A/M$ for all $M$. Therefore $S_{2n}$
is an identity for $B=\Pi_M A/M$, where $M$ ranges over all maximal ideals.
But by the above remarks, $A$ has a natural embedding into  $B$ and hence any identity of $B$ is also an identity of $A$. \qed

\medskip

Let $M$ be a maximal two-sided ideal of $A$. Then $V_A=A/M$ is a simple P.I.
algebra and hence of the form $M_n(D)$, where $D$ is a division ring which is
finite dimensional over its center $C$. Moreover, $D=\End V_A$ (\cite[13.3]{mr}).

If $H$ is a maximal commutative subfield of $D$ then $H=k$ since by
\cite[Theorem 13.10.3 (the proof)]{mr} it is finite dimensional over
$k$ and the latter is algebraically closed. Hence (\cite[Lemma
13.3.4]{mr}), $A/M\cong M_n(k)$ for some $n$.

Thus we conclude

\begin{Prop}\label{piprop} 
 $\deg A=\sup_M\dim_kS$, where $S=A/M$ runs through all
irreducible $A$-modules.
\end{Prop}

\medskip

\begin{Rem}\label{rem-deg}
  The Goldie quotient ring $Q(A)$ of $A$ can be obtained by inverting
  the non-zero {\em central} elements of $A$ \cite[Corollary
  13.6.7]{mr}. Thus, $A$ and $Q(A)$ have the same
  P.I. degree. Therefore, any ring $B$ between $A$ and $Q(A)$ has the
  same P.I. degree as $A$ and in case $B$ is affine over $k$, $A$ and
  $B$ have the same degree. 
\end{Rem}

From \cite[Corollary 13.3.5]{mr} we now get

 \begin{Prop}\label{5}
 $\pid A=(\dim_{Q(Z)}Q(A))^{\frac12}$, where $Q(Z)$
denotes the quotient field of the center of $A$.
\end{Prop}

\medskip

As noted in \cite{jo1} we have

\begin{Prop}\label{pro-reg}
Let $\{a_1,\dots,a_k\}$ be a finite set of regular elements of
$A$. There exists an irreducible representation $\rho$ of $A$ of
maximal degree in which $\rho(a_1),\dots,\rho(a_k)$ are units.
\end{Prop}

\pf Let $B=A[a_1^{-1},\dots,a_k^{-1}]$. Then $A\subseteq B\subseteq
Q(A)$. By Remark~\ref{rem-deg}  the P.I. degrees of $A$ and $B$ are equal.  

Let $\rho$ be an irreducible representation of $B$ of maximal dimension
on a finite-dimensional  vector space $V$ over $k$. Let $\rho^\prime$
denote the restriction of $\rho$ to $A$. Then
$\rho^\prime(a_i)=\rho(a_i)$ for all $i=1,\dots,k$. Moreover, each
$\rho(a_i)$ is a regular linear map, hence by the
Cayley-Hamilton Theorem
its inverse is in  $\rho^\prime(A)$. Thus $\rho^\prime(A)=\rho(B)$ and 
hence $\rho^\prime$ is irreducible.  \qed

\bigskip

A different approach to finding the degree of certain algebras has been found
by De Concini and Procesi \cite[p. 60 \S7]{cp}. We recall some of their results:

\medskip

Let $J=(h_{i,j})$ be a skewsymmetric $n\times n$ matrix such that
 $\forall i,j: h_{i,j}\in{\mathbb Z}$. Given $J$, the quasipolynomial
 algebra $k_J[x_1,\dots,x_n]$ is the algebra over the field $k$
 generated by $x_1,\dots,x_n$ and with defining relations
 \begin{equation}x_ix_j=q^{h_{i,j}}x_jx_i\qquad i<j.  \end{equation} We
 call $J$ the {\bf defining matrix} of the quasipolynomial algebra. In
 the following, $q\in k$ is always assumed to be a primitive $m$th root of
 unity.

\medskip

De Concini and Procesi proved

 \begin{Thm}\label{deppro}
 $\deg k_J[x_1,\dots,x_n]=\sqrt h$, where $h$ is the cardinality
of the image of the map induced by $J$
 \begin{equation}{\mathbb Z}^n\mapsto\left({\mathbb Z}/m{\mathbb Z}
\right)^n,
 \end{equation}
 defined by 
$w=(w_1,\dots,w_n)\mapsto\overline{Jw}$, where 
$\overline{\phantom{w}}$
denotes taking residue class in each coordinate.

Furthermore $k_J[x_1,\dots,x_n]$ is a free module over its center of
rank $\sqrt h$.
\end{Thm}

 In \cite[7.2 Proposition, p.61]{cp} it was shown that
$k_J[x_1^{\pm1},\dots,x_n^{\pm1}]$ is an A\-zu\-maya al\-ge\-bra under
the assumption that $k$ is an algebraically closed field of
characteristic $0$ and  arbitrary $J$.

Since the algebras $k_J[x_1^{\pm1},\dots, x_n^{\pm1}]$ all are affine prime
P.I.  algebras we can use some of the results from above to
prove that the assumptions on $k$ made by De Concini and Procesi are
superfluous. 

\begin{Prop}\label{azu} Let $k$ be a field and $J$ a skewsymmetric matrix with
integer coefficients. The algebra
$k_J[x_1^{\pm1},\dots,x_n^{\pm1}]$ is an Azumaya algebra.
\end{Prop}

\medskip

\pf We wish to use the Artin-Procesi Theorem
\cite[Theorem 13.7.14]{mr}. Thus we have to show that $S_m$
is an identity for $A$ if and only if $S_m$ is an identity for $A/P$
for all primes of $A$ (\cite[Corollary 13.6.7]{mr}).

$S_m$ is an identity for $A$ or $A/P$ if $S_m$ vanishes on all $m$-tuples
of monomials  in the $x_j^{\pm1}$s, $1\leq j\leq n$.

For any such an $m$-tuple  $(a_1,\dots,a_m)$ we get
 \begin{equation}S_m(a_1,\dots,a_m)=f_g\cdot x_1^{l_1}\cdots
   x_n^{l_n}
 \end{equation}
 with $l_j\in{\mathbb Z},\; f_g\in k$ depending on
 $a_1,\dots,a_m$. Since no $x_1^{l_1}\cdots x_n^{l_n}$ can belong to
 a prime $P$, the claim follows. \qed

\bigskip

Many of the algebras considered in the following are iterated Ore
extension. We therefore recall some results on the degree of iterated Ore
extensions or skew polynomial algebras.

In \cite[Theorem 1]{jo1} and in \cite[p. 59, Theorem]{cp} it was
proved that if $R$ is an affine prime algebra over a field of
characteristic $0$, then $\deg R[\theta;\alpha,\delta]=\deg
R[\theta;\alpha]$ provided $\deg R[\theta;\alpha,\delta]$ is finite.
Here, $\alpha$ is $k$-automorphism of $R$ and $\delta$ an
$\alpha$-derivation on $R$.

Combining this result with Proposition~\ref{5}, the degrees of the so called
Dipper-Donkin algebras $D_q(n)$ and the quantized (``official'') matrix algebras $M_q(n)$ were found by Jakobsen and Zhang \cite{jaz1,jaz2}, (for the quantum parameter $q$ a primitive $m$th root of unity).

For later purposes we need a few more results concerning skew
polynomial algebras (implicitly in \cite{jo1}).

We consider an affine prime P.I. algebra $R$ and a skew polynomial
algebra \begin{equation}A=R[\theta;\alpha,\delta],\end{equation} where
$\alpha$ is a $k$-automorphisms of $R$ and $\delta$ an
$\alpha$-derivation.

\medskip

 We assume $A$ is a P.I. algebra and $\alpha$ has finite order.

Notice that for any regular element $r\in R$ (resp. $r\in A$)
$R[r^{-1}]$ (resp. $A[r^{-1}]$) is a   $k$-affine prime P.I.
subalgebra of the Goldie quotient ring of $R$ (resp. $A$) hence by
the previous results has the same degree as $R$ (resp. $A$).

The next lemma follows easily, and is very well known.

\begin{Lem}
\begin{equation}\pid A=\pid A[t].
\end{equation}
\end{Lem}  

Combining these results for $A=R[\theta;\alpha,\delta]$ we obtain
\begin{Lem}\label{cid}
In case there exist regular $r,t$ in $R$ and $s$ in $R$ such that
$r\theta+st^{-1}$ commutes with all elements of $A$ then
\begin{equation}\deg R[\theta;\alpha,\delta]=\deg R.
\end{equation}
\end{Lem}

\medskip

In \cite[Section 4]{jo1} it was proved that such $r,s$ and $t$
exist when $\alpha$ induces the identity on $Z(R)$, the center of
$R$. 

\begin{Rem}
In the present article this result will only be used in situations in
which there exist a regular element $r\in R$ and an element $s\in R$ such 
that $r\Theta +s=z$ is central in $A$. Observe that $r$ is regular in
$A$. In an irreducible representation $\rho$ of $A$ of maximal degree, we may 
by Proposition~\ref{pro-reg} assume that $r$ is invertible. Then
$\rho(\Theta)=\rho(r)^{-1}\rho(z-s)$ and hence $\rho$ remains
irreducible when restricted to $R$. Thus, $\deg A=\deg R$.
\end{Rem}

It was also shown in \cite[the proof of Theorem 3.1]{jo1} that in case
$\alpha$ is not the identity on $Z(R)$, then there exists a finitely
generated multiplicatively closed $\alpha$-invariant set $T$ of
central elements of $R$ such that
 \begin{equation}R[T^{-1}][\theta;\alpha,\delta]\simeq
   R[T^{-1}][\theta';\alpha], 
 \end{equation}
 where $\theta'=\theta-a$ for some $a\in Z(R[T^{-1}])$, and such that
$R[T^{-1}]$ is $k$-affine.

\medskip

In case there exists a subalgebra $Z_0$ of the center $Z$ of $R$ such
that i) $Z$ is a finite $Z_0$ module, ii) $\delta(Z_0)=0$, and iii)
$\alpha_{\mid Z_0}=1_{Z_0}$, De Concini and Procesi proved
(\cite[Theorem p. 58]{cp}) that
\begin{equation}
\deg R[\Theta,\alpha,\delta]=(\deg R)\cdot k,
\end{equation}
where $k$ is the order of $\alpha$'s restriction to $Z(R)$.

\smallskip

[By the methods of \cite{jo1} and \cite{jo2} one can in fact show that the
special assumptions on $R,\alpha,\delta$ are superfluous. One just
needs that $R[\Theta,\alpha,\delta]$ is a prime P.I. algebra.]

\smallskip

In particular, we get, provided $R[\Theta,\alpha,\delta]$ is a P.I. algebra,

\begin{Prop}\label{coa}
  Let $R$ be a prime P.I. algebra and $\alpha$ an automorphism of $R$
  of order $k$. If there exists an element $c\in Z(R)$ such that the
  $\alpha$ orbit of $c$ has order $k$, then
\[\deg (R[\Theta,\alpha,\delta])=(\deg R)\cdot k.
\]
\end{Prop}

\medskip

\section{The quantized function algebra $M_q(n)$}\label{stand}

The ``standard'' quantized function algebra $M_q(n)$ of $n\times n$
matrices is the quadratic algebra generated by $n^2$ elements
$Z_{i,j},\;i,j=1,\dots,n$ and with defining relations
  \begin{eqnarray}\label{def-rel}Z_{i,j}Z_{i,k}&=&qZ_{i,k}Z_{i,j} \text{ if } j<k,\\\nonumber
Z_{i,j}Z_{k,j}&=&qZ_{k,j}Z_{i,j} \text{ if }i<k,\\\nonumber
Z_{i,j}Z_{s,t}&=&Z_{s,t}Z_{i,j} \text{ if }i<s,t<j,\\\nonumber
Z_{i,j}Z_{s,t}&=&Z_{s,t}Z_{i,j}+(q-q^{-1})Z_{i,t}Z_{s,j} \text{ if }
i<s, j<t,
\end{eqnarray}
for $i,j,k,s,t=1,2,\cdots,n$.  

It is well known that the monomials $Z^A=Z_{1,1}^{a_{1,1}}\cdots  Z_{1,n}^{a_{1,n}}\cdots
Z_{2,1}^{a_{2,1}}\cdots  Z_{2,n}^{a_{2,n}}\cdots
Z_{n,1}^{a_{n,1}}\cdots  Z_{n,n}^{a_{n,n}}$ for
$A=\{a_{i,j}\}_{i,j=1,1}^{n,n}\in Mat(n^2,{\mathbb N}_0)$ form a PBW-type basis
of $M_q(n)$ for any
$q\neq0$. The quantum determinant 
\begin{equation} {\det}_q=\Sigma_{\sigma\in
    S_n}(-q)^{l(\sigma)}Z_{1,\sigma(1)}Z_{2,\sigma(2)} \cdots
  Z_{n,\sigma(n)}.\end{equation}is central for any $0\neq
q\in k$.

Viewing $M_q(n)$ as an iterated skew polynomial algebra (c.f. below),
it follows that {\it the associated quasipolynomial algebra}
$\overline{M_q(n)}$ is given in terms of the same generators, but with
defining relations
\begin{eqnarray}Z_{i,j}Z_{i,k}&=&qZ_{i,k}Z_{i,j} \text{ if } j<k,\\\nonumber
Z_{i,j}Z_{k,j}&=&qZ_{k,j}Z_{i,j} \text{ if }i<k,\\\nonumber
Z_{i,j}Z_{s,t}&=&Z_{s,t}Z_{i,j} \text{ if }i<s,t<j,\\\nonumber
Z_{i,j}Z_{s,t}&=&Z_{s,t}Z_{i,j} \text{ if } i<s, j<t,
\end{eqnarray}
for $i,j,k,s,t=1,2,\cdots,n$.

\medskip

Later on, we shall encounter a number of subalgebras $B$ of
$M_q(n)$. For each of these, analogously to the above, the associated
quasipolynomial algebra $\overline{B}$ is the algebra with the same
generators but where in the defining relations, all terms of the form
$(q-q^{-1})Z_{i,t}Z_{s,j}$ have been dropped.

\medskip

The degree of $M_q(n)$ was found in \cite{jaz1} to be $m^{n(n-1)/2}$
 for $q$ an $m$th root of unity, $m$ odd. The approach there utilized
 the result of Procesi and De Concini (\cite[Theorem p. 59]{cp}) according to which, as
 a special case, $\deg M_q(n)=\deg\overline{M_q(n)}$.  We will reprove
 this result by utilizing the results of Section~\ref{pi}.

Before doing so let us introduce some notation, which will be used
also later in this paper.

We view $M_q(n-1)$ as  the $k$-algebra generated by the elements
$Z_{i,j}$, $1\leq i,j\leq
n-1$ and we will then view $M_q(n)$ as an iterated Ore extension of
$M_q(n-1)$ obtained by adjoining the indeterminates as follows:

\begin{equation} M_q(n)=M_q(n-1)[Z_{n,1};\alpha_{n,1}]\cdots[Z_{n,n-1};
\alpha_{n,n-1}\delta_{n,n-1}] [Z_{1,n};\alpha_{1n}]\cdots[Z_{n,n};
\alpha_{n,n},\delta_n],
 \end{equation}
where $Z_{n,k}Z_{i,j}=\alpha(Z_{i,j})Z_{n,k}+\delta(Z_{i,j})$ for
 $1\leq i,j<n$ { or } $i=n,j<k$, and where
$Z_{k,n}Z_{i,j}=\alpha(Z_{i,j})Z_{k,n}+\delta(Z_{i,j})$ for $1\leq i,j<n$ or
$j=n,i<k$.

Changing slightly the notation from Parshall and Wang \cite[Section
4]{pw} we define $I=\{n-i+1,\dots,n\}$ and $J=\{1,\dots,i\}$ and let
$\widetilde\theta_{n+1-i}=D(I,J)$ (the quantum determinant based on
the rows in $I$ and columns in $J$) and $\theta_{i+1}=A(I,J)$. Later on,
we shall introduce some more general elements with these names but
letting $r=n$ in Figure~\ref{f1}, the elements
$\theta_k,\widetilde\theta_t$ in that figure are precisely what have been
defined here (with $k=i+1$ and $t=n+1-i$).

In \cite{jaz1} the following notion was introduced

\begin{Def}\label{covdef} An element $x\in M_q(n)$ is called covariant if for any
$Z_{i,j}$ there exists an integer $n_{i,j}$ such that
\begin{equation}xZ_{i,j}=q^{n_{i,j}}Z_{i,j}x.\end{equation}
Clearly, $Z_{1,n}$ and $Z_{n,1}$ are covariant.
\end{Def}

It was then shown that  the elements $\widetilde\theta_{t}$
and $\theta_{k}$ are covariant. Utilizing this, the following elements were also found to
belong to the center: 
 \begin{equation}c_{i+1}=\widetilde\theta_{n+1-i}\theta_{i+1}^{m-1}\mbox{ and
 }d_{i+1}=\widetilde\theta_{n+1-i}^{m-1}\theta_{i+1} \end{equation} 
 for $1\leq i\leq n-1$.
 
For the convenience of the reader we list the covariance properties of
the elements $\theta_{i+1},\widetilde\theta_{n+1-i}$ (from which the
centrality of $c_{i+1}$ and $d_{i+1}$ also is obvious). Let
$I=\{n-i+1,\dots,n\}$ and $J=\{1,\dots,i\}$ as previously:
\begin{eqnarray}
  \label{eq:covrel}
  Z_{a,b}\theta_{i+1}=q\theta_{i+1} Z_{a,b}&\textrm{ and
    }&Z_{a,b}\widetilde\theta_{n+1-i}=
  q\widetilde\theta_{n+1-i}Z_{a,b}\textrm{ for }a\notin I, b\in J, \\
Z_{a,b}\theta_{i+1}=      q^{-1}\theta_{i+1} Z_{a,b}&\textrm{ and
    }&Z_{a,b}\widetilde\theta_{n+1-i}=
  q^{-1}\widetilde\theta_{n+1-i}Z_{a,b}\textrm{ for }a\in I, b\notin
  J, \\  
Z_{a,b}\theta_{i+1}= \theta_{i+1} Z_{a,b}&  \textrm{ and  }&Z_{a,b}\widetilde\theta_{n+1-i}=
  \widetilde\theta_{n+1-i}Z_{a,b}\textrm{ in all other cases}. 
\end{eqnarray}

\medskip

It is well-known that $M_q(n)$ (being an iterated Ore extension) is a domain, thus the $r,t$ in
 Lemma~\ref{cid} are automatically regular if non-zero.

\label{detsec}

\begin{Thm}\label{detdeg} $\deg M_q(n)=m^{n(n-1)/2}$, where $m$ is an odd integer
and $q$ is a primitive $m$th root of unity.
\end{Thm}

\pf We use induction on $n$. Since the formula clearly holds for $n=1$, it
suffices to prove \begin{equation}\deg M_q(n)=m^{n-1}\deg M_q(n-1).
\end{equation} First notice that
\begin{equation}{\det}_q=rZ_{n,n}+s, \end{equation} where in fact $r$
up to a sign is $\det_q$ for $M_q(n-1)$,  and where $s$, when expanded in
the PBW basis, does not contain $Z_{n,n}$ either.

By Lemma~\ref{cid} we see that $\deg M_q(n)$ is the same as the degree of
the algebra where $Z_{n,n}$ is excluded.

The same argument works for $Z_{n-i,n}$. One just has to replace $\det_q$ by $d_{i}=\widetilde\theta_{n+1-i}^{m-1}\theta_{i+1}$ in the procedure. Thus, 
 \begin{equation}\deg M_q(n)=\deg M_q(n-1)[Z_{n,1};
\alpha_{n,1}]\cdots[Z_{n,n-1};\alpha_{n,n-1},\delta_{n,n-1}].
\end{equation}

Let $R_j$ be the algebra obtained by adjoining $Z_{n,1},\dots, Z_{n,j}$
to $M_q(n-1)$. Let $\underline{\widetilde\theta}_{i}$ and $\underline{\theta}_i$ be the
quantities in $M_q(n-1)$ analogous to $\widetilde\theta_i$ and $\theta_i$, Then
notice that \begin{equation}\underline{
c}_{j+1}=\underline{\widetilde\theta}_{n+1-j}\underline{\theta}_{j+1}^{m-1} (\underline{{\det}}_q)^{m-1}
\end{equation}
is a central element in $R_j$ and 
 \begin{equation}\alpha_{n,j+1}(\underline{c}_{j+1})=q^2\underline{c}_{j+1}.
 \end{equation}

Therefore by Proposition~\ref{coa} and because $m$ is odd, $\deg
R_{j+1}=m\deg R_j$.

In the case of $R_1$ one may just use $\underline{\det}_q$ as $c_0$. Since
$\alpha_{n,1}(c_0)=q^{-1}c_0$, there is a similar conclusion. The proof is thus completed. \qed

\medskip

\section{The case of $n\times r$}\label{nr}

We consider here the quadratic algebra $M_q(n,r)$ consisting of $n\times r$
matrices ($r\leq n$). We determine the degree and the center in case
$m$ is ``good'' and also obtain  insight, in some cases even full,  into the cases where $m$
is not so ``good''. 

\subsection{\underline{Central elements}}

We assume that $r=x\cdot s$ and $n=y\cdot s$ with $x\cdot y$ odd and $s$ as
big as possible. We display $s$ central elements of $M_q(n,r)$. (For a 
hint of how these were discovered, see the proof of
Proposition~\ref{bl2} below. Also recall from \cite{jaz1} that the
centres of $M_q(n,r)$ and $\overline{M_q(n,r)}$ are in bijective
correspondence via the leading term.)

We begin by defining elements $\Psi_t\in M_q(n,r)$ for $t=1-r,\dots,
n$. First set $\Psi_{1-r}=1$.  $\Psi_{2-j}$ is the quantum $(r-j+1)\times
(r-j+1)$ minor involving the rows $1,2,\dots, r-j+1$ and columns
$j,j+1,\dots,r$ (for $j=2,\dots,r$). Then consider $i=1,\dots, n-r+1$ and let
$\Psi_i$ be the quantum $r\times r$ determinant involving the rows
$i,i+1,\dots,i+r-1$. Finally, $\Psi_{n-k+1}$ is the
quantum $k\times k$ determinant involving columns $1,\dots,k$ and rows
$n-k+1,n-k+2,\dots,n$ for $k=1,\dots, r-1$.

\begin{Lem}
For $a=1,\dots,s$ the elements $Z_a$ defined by
\begin{equation}\label{aeq}
Z_a:=\prod_{\ell=-x}^{y-1} (\Psi_{a+\ell\cdot s})^{(-1)^\ell}
\end{equation}
are central.
\end{Lem}

\pf Let us consider the case $a=1$. We may then view our $n\times r$ matrix as
being built up of $y\cdot x$ blocks $B_{i,j}$ of size $s$, block $B_{1,1}$
consisting of rows and columns $1,\dots,s$, block $B_{2,1}$ consisting of rows
$s+1,\dots 2s$ and columns $1,\dots,s$, etc. Let us now look at some
$X_{a,b}\in B_{i,j}$.  Due to the covariance of the various determinants it is
possible to see that the commutativity of $Z_1$ with $X_{a,b}$ is equivalent
to picking up a factor of $q^{\pm1}$ for each each block $B_{\alpha,j}$ and
each $B_{i,\beta}$ and that indeed the whole computation may be viewed as the
computation for commutativity of the analogous expression computed in
$\overline{M_q(y,x)}$.  Here it is a matter of investigating the matrix
$B=\{b_{i,j}\}_{i,j=1}^{y,x}$ given by $b_{i,j}=(-1)^{i+j}$ and checking that
${\overline{Z}}_1^B$ is central. But since $x$ and $y$ are odd, this is
straightforward. Indeed, the computation is reduced to ascertaining that if
$w$ is odd and $1\leq i\leq w$ then $(w-i)-(i-1)=0$ in ${\mathbb Z}_2$.

\begin{figure}
\epsffile{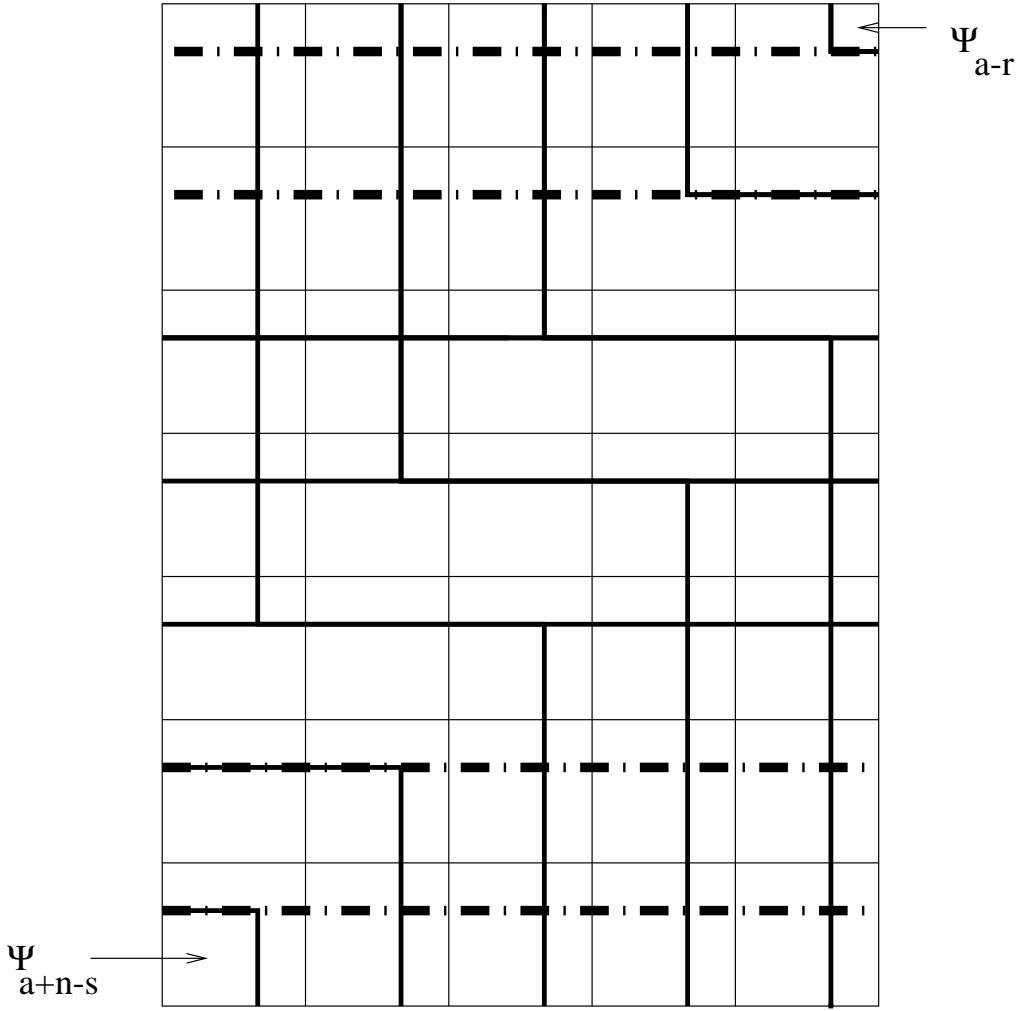}
\caption{\label{fig6}The factors $\Psi_i$ for $a\geq2$. Beginning at the 
  upper right hand corner, the first is an $(a-1)\times (a-1)$ matrix,
  then they grow larger in steps of size $s$, reach full size as
  indicated by the dashed lines, and then decrease again until
  reaching the bottom left corner which is a minor of size
  $(s-a+1)\times (s-a+1)$. The thin lines represent blocks of size
  $s\times s$. At $a=1$ the factor $\Psi_{1-r}$ deteriorates to a 1
  and the lines defining the factors run into those defining the
  blocks.  The reader is now invited to take a stroll in the figure:
  Place yourself in the position of some $Z_{i,j}$, say $Z_{n,1}$, and
  count covariance-$q$'s from each factor. Convince yourself that the
  exponents are exactly right for giving a total factor $q^0$. Then
  move to a neighboring element and so forth. When crossing
  a ``boundary'' to a new factor, there will precisely be one new $q$
  and one new $q^{-1}$.}
\end{figure}

The case $a\geq2$ is similar though slightly more complicated. One may
view the relevant diagram over the minors as being obtained by a
translation along the line from the upper right corner to the lower
left. Only near the boundaries does one need to check carefully that
the signs match up. The details are left to the reader,
c.f. Figure~\ref{fig6}. \qed

\medskip

Finally observe
that one has the following result,  the first part of which is  just as in \cite{jaz1}:

\begin{Lem}\label{5.222} i) Let $m$ be even and $q$ a primitive $m$th root of unity. Any
  element of the form
\begin{equation}\label{521}
Z_{k,i}^{\frac{m}2}   Z_{k,j}^{\frac{m}2}   Z_{\ell,i}^{\frac{m}2}
 Z_{\ell,j}^{\frac{m}2},
\end{equation} 
where $1\leq k<\ell\leq n$ and $1\leq i<j\leq r$, is in the center of
$\overline{M_q(n,r)}$. 

ii) If $n+r$ is even then 
\begin{eqnarray}\label{522}
 \left(Z_{n,1}\dots Z_{2,1}Z_{1,1}Z_{1,2}\dots
   Z_{1,r}\right)^{\frac{m}2}&&\textrm{ is central if }n,r \textrm{
   are odd, and}\\\nonumber
 \left(Z_{n,1}\dots Z_{2,1}Z_{1,2}\dots
   Z_{1,r}\right)^{\frac{m}2}&&\textrm{ is central if }n,r \textrm{
   are even.}
\end{eqnarray}

\end{Lem}

\begin{Rem}
The analogous elements of $M_q(n,r)$ (i.e. the same expressions but
where the generators are in that algebra) are  also central in this case.
\end{Rem}

\begin{Rem}
It follows by an argument similar to the one in
\cite[Theorem~6.2]{jaz1} that in case $m$ is ``good''
(c.f. (\ref{goodlabel}) and Proposition~\ref{gcd} below) then the elements
$Z_a;a=1,\dots,s$ together with the elements $Z_{i,j}^m; 1\leq
i\leq n; 1\leq j\leq r$ generate the center.
\end{Rem}

\medskip

\subsection{\underline{Degree and block diagonal form}}

\medskip

Following \cite[p. 469-470]{jaz1}, the defining matrix $J$ of the associated quasipolynomial algebra
$\overline{M_q(n,r)}$, with respect to a natural basis $\{E_{i,j}\}$, is in fact the matrix of the map
\begin{equation}\label{44}
A\stackrel{J}{\rightarrow} H_n A-AH_r,
\end{equation}
where $H_k=\sum_{1\leq j<i\leq k}(E_{i,j}-E_{j,i})$ for $k=n,r$. We start by
computing the rank of this map or, equivalently, the dimension of the
kernel. Let
\begin{equation}
c_{n,r}:=\corank(J).
\end{equation}
Thus, if $m$ is ``good'', e.g. a large prime,
\begin{equation}\label{goodlabel}
\deg(M_q(n,r))=m^{\frac12(nr-\corank)}=m^{\frac12(nr-c_{n,r})}.
\end{equation}
Indeed, $\deg(M_q(n,r))=\deg(\overline{M_q(n,r)})$ by \cite[Theorem
p. 59]{cp}, and the latter is given by Theorem~\ref{deppro}. The
cardinality of $J$ is clearly the same as for any block diagonal form
of $J$. Now consider a $2\times 2$ matrix
$\left(\begin{array}{cc}0&r\\-r&0\end{array}\right)$ for some
$r\in\mathbb N$.  Considered as a map from ${\mathbb Z}^2$ to $({\mathbb Z}/m{\mathbb
  Z})^2$, the image clearly has cardinality $(m/\textrm{g.c.d.}(m,r))^2$, where
$\textrm{g.c.d.}(m,r)$ denotes the greatest common divisor of $m$ and
$r$. This means that if e.g. $m$ is a prime bigger than all
non-zero elements in a block form of $J$ then each non-trivial block
contributes with a factor $m$ to the total degree. But clearly, there
are $\frac12\cdot \rank J$ such blocks.

\medskip

Recall from \cite[p. 470]{jaz1} that 
$H_k=S_k+\dots+S_k^{k-1}=\frac{1+S_k}{1-S_k}$ where
$S_k=-E_{1,k}+\sum_{i=2}^kE_{i,i-1}$. In particular, $S_k^k=-1$. Now
observe that
\begin{equation}\label{m-eq}
  H_n A-AH_r=M^\prime\Leftrightarrow 2\left(S_nA-AS_r\right)=(1-S_n)M^\prime (1-S_r),
\end{equation}
where $1-S_k$ is invertible, indeed,
$(1-S_k)^{-1}=\frac12(1+S_k+\dots+S_k^{k-1})$.

This observation will also be used later, but it follows immediately
that the kernel is given by those $n\times r$ matrices $A$ for which
\begin{equation}\label{rela}
  S_nA=AS_r.
\end{equation}
If $A=\sum_{i,j}a_{i,j}E_{i,j}$ a straightforward computation gives that
(\ref{rela}) is equivalent to 
\begin{eqnarray}\label{newform1}
  \forall i=2,\dots,n,\forall j=1,\dots,r-1:
  a_{i-1,j}&=&a_{i,j+1},\\\nonumber
\forall i=2,\dots,n: a_{i-1,r}&=&-a_{i,1},\\\nonumber
\forall j=1,\dots,r-1: a_{1,j+1}&=&-a_{n,j},\textrm{ and }\\\nonumber
a_{n,r}&=&a_{1,1}.
\end{eqnarray}\label{ext}
If we define, for all $\alpha,\gamma\in\mathbb Z$,   for
$\beta=0,\dots,n-1$, and for  $\delta=0,\dots,r-1$
\begin{equation}\label{newform2}
  a_{\alpha n+\beta,\gamma r+\delta}=(-1)^{\alpha+\gamma}a_{\beta,\delta},
\end{equation}
then (\ref{rela}) is equivalent to
\begin{equation}\label{trel}
 \forall \beta,\delta,t:\qquad\
 a_{\beta+t,\delta+t}=a_{\beta,\delta}.
\end{equation}

\medskip

\begin{Prop}\label{gcd}
Let $s=\textrm{g.c.d.}(n,r)$. Specifically, let $n=x\cdot s$ and $r=y\cdot
s$. Then $J$ is non-invertible if and only both $x$ and $y$
are odd. In this case,  $$c_{n,r}=\corank J=s.$$
\end{Prop}

\pf Observe that by definition, $x$ and $y$ cannot both be even. Now,
according to (\ref{newform2})  and (\ref{trel}) 
\begin{equation}
 a_{y\cdot n,x\cdot r} =(-1)^{x+y}a_{0,0}=(-1)^{x+y}a_{1,1}=a_{x\cdot
   y\cdot s, x\cdot y\cdot s}=a_{1,1}.
\end{equation}
Thus, if a solution to (\ref{rela}) is to exist with $a_{1,1}\neq 0$
then $(-1)^{x+y}=1$. In this case, a solution is given by $\forall
t:a_{t,t}=a_{1,1}=1$ and all other $a_{i,j}=0$. More generally, a
non-zero solution exists if and only if $(-1)^{x+y}=1$. In this case
there are $s$ independent solutions given by
\begin{equation}
 \forall t: a_{t+i,t}=a_{1+i,i}\textrm{ for }i=0,1,\dots,s-1 \textrm{
   and all other } a_{i,j}=0. 
\end{equation}
The claim now follows. \qed

\medskip

We wish to obtain more precise information about the blocks of a
diagonal form of the associated matrix of the quasipolynomial algebra
$\overline{M_q(n,r)}$.

\medskip

\begin{Prop}\label{bl1}
In the case where $q=-1$ there is an irreducible module of dimension
$2^{d_0}$ where 
\[d_0={\left[\frac{n+r-1}2\right]}.\]
\end{Prop}

\pf In this case the term $(q-q^{-1})$ disappears and so
$\overline{M_{-1}(n,r)}=M_{-1}(n,r)$. Hence, by covariance, in an irreducible
module, $Z_{i,j}$ is either zero or invertible. Consider the subalgebra 
\begin{equation}
{\s S}_{n,r}=\left<Z_{1,j},Z_{i,1}\mid 1\leq i\leq n\textrm{ and } 1\leq j\leq r\right>.
\end{equation}
By Proposition~\ref{pro-reg} there is
an irreducible representation of maximal dimension of this algebra in
which all the generators are invertible. Given such an irreducible
module, the recipe
\begin{equation}
Z_{i,j}=c_{i,j}\cdot Z_{1,1}^{-1}Z_{i,1}Z_{1,j},
\end{equation}
where $2\leq i\leq n$, $2\leq j\leq r$, and the $c_{i,j}$s are
arbitrary constants, define an irreducible representation of
$\overline{M_q(n,r)}$. The proof is thus completed if we can establish
that there are precisely $\left[\frac{n+r-1}2\right]$ non-trivial
blocks in the block diagonal form of the associated matrix of ${\s
  S}_{n,r}$ (and hence that the degree is
$m^{\left[\frac{n+r-1}2\right]}$). For this purpose, let $x_j=Z_{1,j}$
for $j=1,\dots,r$ and let $y_j=Z_{i,1}$ for $i=2,\dots,n$. Upon the
replacements $x_j\mapsto x_1x_2x_j$ ($j=3,\dots,r$), $y_2\mapsto
y_2x_2$, and $y_i\mapsto y_2y_i$ ($i=3,\dots,n$), the pair $x_1,x_2$
decouples completely leaving us with an algebra which is isomorphic to
${\s S}_{n-1,r-1}$. It is well-known (and elementary) to see that
there are $\left[\frac{x}2\right]$ non-trivial blocks in the block
diagonal form corresponding to ${\s S}_{x,1}$. The result follows
directly from these observations. \qed

\medskip

Actually, we did not use anything about the algebra except that it was contained in a box of size $n\times r$, hence we get the following corollary to the proof:

\begin{Cor}
Let $S$ be a subalgebra of $M_q(n,r)$ such that $\forall i=1,\dots,n:
Z_{i,1}\in S$ and such that $\forall i=j\dots,r: Z_{1,j}\in S$. Then,
in case $q=-1$ there is an irreducible module of dimension $2^{d_0}$
where
\[d_0={\left[\frac{n+r-1}2\right]}.\]
\end{Cor}

\medskip

\begin{Rem}
  For general $q$ we get a result similar to Proposition~\ref{bl1}.
  Specifically, given a representation of $S_{n,r}$ in which the
  generators are denoted $\overline{Z}_{i,j}$ and in which
  $\overline{Z}_{1,1}$ is invertible, the recipe
\begin{eqnarray*}
    Z_{i,j}&=&q\overline{Z}_{1,1}^{-1}\overline{Z}_{i,1}\overline{Z}_{1,j}\textrm{ if }i,j>1\\
Z_{i,j}&=&\overline{Z}_{i,j}\textrm{ else }
  \end{eqnarray*}
defines a representation of $M_q(n)$ as can be seen by a
straightforward but tedious computation.
\end{Rem}

\medskip


\subsection{\underline{The general form of the center and the blocks}}

We wish to take a closer look  at the degree and center in the case where $m$ is not necessarily a 
prime and where $n,r$ are arbitrary. For this purpose we need more
information about the pertinent block diagonal form.

\begin{Prop}\label{bl2}
The non-trivial blocks in a block diagonal form of the defining matrix
 $J$ of $\overline{M_q(n,r)}$ are either of the form
 $\left(\begin{array}{cc}0&1\\-1&0\end{array}\right)$,
 $\left(\begin{array}{cc}0&2\\-2&0\end{array}\right)$, or of the form
 $\left(\begin{array}{cc}0&4\\-4&0\end{array}\right)$.
\end{Prop}

\pf We will begin by studying the center of $\overline{M_q(n,r)}$ at a
primitive $m$th root of unity. If $A$ is an $n\times r$ integer
matrix, the condition for a monomial
$u^A=Z_{1,1}^{a_{1,1}}Z_{1,2}^{a_{1,2}}\dots Z_{n,r}^{a_{n,r}}$ to be in
      the center is precisely (in the notation of (\ref{44})) that
\begin{equation}\label{mod-eq}
  H_n A-AH_r=0\mod m. 
\end{equation}
Returning to (\ref{m-eq}), it follows that 
\begin{equation}\label{symeq}
 S_nA-AS_r={\frac12}\cdot M,
\end{equation}
where $M=(1-S_n)M^\prime(1-S_r)$ is an integer matrix whose entries
all are multiples of $m$. However, basically due to the $\frac12$ in
$(1-S_k)^{-1}$, not all such matrices $M$ need define a solution $A$
to (\ref{mod-eq}). Returning now to the equations (\ref{newform1} -
\ref{trel}), these remain valid when reinterpreted as equations {\em
  modulo $\frac{m}2$}. In case $(-1)^{x+y}=1$, with $x,y$ as in
Proposition~\ref{gcd}, we just get the old solutions possibly with
some elements of the form (\ref{521}) or (\ref{522}) superimposed. But
the case $(-1)^{x+y}=-1$ now is possible since we just need
$2a_{1,1}=0\mod \frac{m}2$. Thus, it should be proportional
to $\frac{m}4$. In all cases it follows that the entries of
$A$ are integer multiples of $\frac{m}4$ and hence that the element
$u^A$ satisfies that its fourth power is in the central subalgebra
generated by the elements $Z_{i,j}^m$. But suppose that there is a
block of the form $\left(\begin{array}{cc}0&s\\-s&0\end{array}\right)$
with $s\neq1,2,4$. Then there are monomials $u^A,u^B$ such that
$u^Au^B=q^su^Bu^A$ and such that $u^A$ commutes with all other
generators. But then the element $(u^A)^{\frac{m}s}$ is central for
any $m$ which is a multiple of $s$ and this is a contradiction since
$(u^A)^{\frac{4m}s}$ will not be in the above mentioned central
subalgebra.\qed

\medskip

\begin{Rem}\label{remmm}
  It follows from (\ref{symeq}) that if $u^A$ is central, then so is
  $u^B$ for any $B=S_n^i A S_r^j$ with $i,j\in{\mathbb Z}$. This
  symmetry can be used to construct new solutions from given ones,
  c.f. below. 
\end{Rem}

\medskip

\begin{Prop}\label{4.6}
  The non-trivial blocks in a block diagonal form of the defining
  matrix $J$ of $\overline{M_q(n,r)}$ are: $d_0$ matrices of the form
  $\left(\begin{array}{cc}0&1\\-1&0\end{array}\right)$ and \newline
  $\max\{0, \frac{nr-c_{n,r}}{2} - d_0\}$ matrices of the form
  $\left(\begin{array}{cc}0&2\\-2&0\end{array}\right)$ or 
  $\left(\begin{array}{cc}0&4\\-4&0\end{array}\right)$.
\end{Prop}

\pf This follows immediately from Proposition~\ref{azu}, Proposition~\ref{bl1}, and
Proposition~\ref{bl2}. \qed
 
\medskip

We finish with some remarks about the occurrence of ``$4$''s: Suppose
that $r,n$ are relatively prime. Following the reasoning just before
Remark~\ref{remmm}, if we are to have a genuine solution to
(\ref{symeq}) involving $\frac{m}{4}$ then by (\ref{trel}) we must
have $a_{\beta ,\delta }=$ $\frac{m}{4}\mod $ $\frac{m}{2} $ for all
$\beta ,\delta $. We may assume that $a_{\beta ,\delta }=$ $\frac{
  m}{4}$ for all $\beta =2,\ldots ,r$ and all $\delta =2,\ldots ,n$.
This follows since if some $a_{\beta ,\delta }=$ $\frac{3m}{4}$ then,
by Lemma ~\ref{5.222} the $\frac{m}{2}\,$is part of a central element
involving $a _{1,1},a_{1,\delta },a_{\beta ,1}$, and $a_{\beta ,\delta
  }$ and may thus be discarded. We furthermore assume that
$a_{1,j}=\frac{m}{4} +\alpha _{j}\frac{m}{4}$ for $j=1,\ldots ,n$ and
$a_{i,1}=\frac{m}{4}+\beta _{i}\frac{m}{4}$ for $i=2,\ldots ,r$, where
each $\beta _{i}\,$and $\alpha _{j}$ is $0$ or $2$ modulo $4$.

Consider first a pair of indices $(i,j)\,$with $i,j>1$. Then the condition
for $u^{A}\,$to commute with $Z_{i,j}\,$in the quasipolynomial algebra is 
\[
(r-i)+(n-j)-(i-1)-(j-1)-\alpha _{j}-\beta _{i}=0\mod 4.
\]
By subtracting consecutive terms it follows that 
\begin{eqnarray*}
\alpha _{j} &=&2j+c\text{ for }j=2,\ldots ,n\text{, } \\
\beta _{i} &=&2i+d\text{ for }i=2,\ldots ,r\text{, and} \\
\alpha _{1} &=&f\text{.}
\end{eqnarray*}
Also observe that 
\[
n+r\,\text{ must be even}
\]
and 
\[
n+r+2+c+d=0\mod 4\text{.}
\]
At a point $(1,j)$ with $j>1$ we get, utilizing the parity properties, 
\[
n+r+nc+n(n+1)+f=0\mod 4\text{.}
\]
Likewise, at $(i,1)$ with $i>1$ we get 
\[
n+r+rd+r(r+1)+f=0\mod 4\text{,}
\]
and, finally, at $(1,1)\,$we get 
\[
n+r-2+(n-1)c+(r-1)d+n(n+1)+r(r+1)=0\mod 4\text{.}
\]

For these equations to have solutions, $r$ and $n$ must have the the
same parity. If both are odd, there are no further restrictions for
these equations to have a solution which also solves (\ref{symeq}). 

Returning to the general situation, if both $n,r$ are even, it turns
out to be a further necessary (and also sufficient) condition for
(\ref{symeq}) to have a genuine $\frac{m}4$ solution that they are
equal modulo $4$. It should now be observed that if there is a general
center, that is, if $J$ is singular, amongst these solutions there
will be some which are actually of no interest, namely those solutions
that correspond to a zero on the right hand side of (\ref{symeq}).

Let us now assume that $r$ is prime.

If $r=2$, $n$ odd does not contribute by the above analysis. If it is
even, it is forced to be of the form $n=4t+2$ which means that the
central elements already have been picked up by the general central
elements.  Suppose then that $r$ is an odd prime. If $r$ does not
divide $n$, is it also forced to be odd and hence $J$ is singular and
hence the solutions we pick up are general central elements raised to
the power $\frac{m}4$.  This takes care of
all cases except $n=zr$ for some positive integer $z$.  If $z$ is odd
we have an $r$-dimensional center: again nothing new.  Finally, if $z$
is even, the previous elements do not give anything.  However, there
are in fact some non-trivial $\frac{m}{4}$-central elements.
Specifically, let $A_{1}$ be the matrix whose non-zero coefficients
$a^{(1)}_{i,j}$ satisfy
\begin{eqnarray*}
 & &a^{(1)}_{1,r}=a_{i,i+jr+e}^{(1)}=\frac{m}{4}\text{ \quad for }\\\nonumber & & i=1,\ldots
r;j=0,\ldots ,z-1, e=0,1\text{,
and }i+jr+e\leq z\cdot r\text{,}
\end{eqnarray*}

and add to that the matrix $A_{2}$ whose whose non-zero coefficients
$a^{(2)}_{i,j}$ satisfy
\[
a_{i,1}^{(2)}=\frac{m}{2}\,=a_{r,1+jr}^{(2)}\text{ for }i=1,\ldots r\text{ and }
j=1,\ldots z-1\text{.}
\]
Then $A_{1}+A_{2}$ defines a central element. Moreover, using the
symmetries of the original equation, we get in fact $(r-1)$ solutions
(which clearly is as much as could be hoped for).

\medskip

We thus have the following partial result

\begin{Prop}
Let $r$ be a prime. Then
$\left(\begin{array}{cc}0&4\\-4&0\end{array}\right)$ occurs in the
block diagonal form  of the defining
  matrix $J$ of $\overline{M_q(n,r)}$ if and only if $r$ is odd and
  $n=z\cdot r$ for some even integer $z$. In this case, there are
  $\frac{r-1}2$ such blocks.
\end{Prop}

\medskip

\section{Quantized minors}\label{qm}


For each $\ell=1,\dots, n$, $I^{\ell}_q$ denotes the ideal generated
by all $\ell\times \ell$ quantum minors. We consider here the
function algebra of rank $r$ matrices. Specifically, let
\begin{equation}
M^{r+1}_q(n)=M_q(n)/I^{r+1}_q.
\end{equation}
For each $t=1,\dots,n$, let $d_t$ denote the $t\times
t$ quantum determinant of the subalgebra generated by the elements
$Z_{i,j}$ with $1\leq i,j\leq t$. Set $d=d_r$.  The natural candidate for quantized
rank $r$ matrices is then 
\begin{equation}\label{sense}
M^{r+1}_q(n)[d^{-1}],   
\end{equation}
where we shall return to the issue of inverting $d$ shortly. 

We wish to compare this algebra to a somewhat more manageable one, namely $A_{n,r}$, where
\begin{Def}The algebra  $A_{n,r}$ is the subalgebra of $M_q(n)$
generated by those $Z_{i,j}$ for which 
$(i,j)\notin\{r+1,\dots,n\}\times\{r+1,\dots,n\}$.
\end{Def}

In \cite{gl}, Goodearl and Lenagan proved that $I^{r+1}_q$ is
completely prime. Moreover, Rigal proved that
$A_{n,1}[d^{-1}]\simeq  M^2_q(n)[d^{-1}]$ (\cite{ri}). We shall prove
below that 
\begin{equation}\label{sj-star}A_{n,r}[d^{-1}]\simeq
  M^{r+1}_q(n)[d^{-1}]
\end{equation} for a general 
$r$ and on the way give a new proof of the former result. 

For  (\ref{sense}) to make sense we first of all need the
following:

\begin{Prop}\label{5.2}
$d=d_r$ is regular in $M_q(n)/I^{r+1}_q$ both in
the case of $q$ generic and the case where $q=\varepsilon$ is a
primitive $m$th root of unity. 
\end{Prop} 

\pf It is proved in \cite[Theorem 2.5]{gl} that $M_q(n)/I^{r+1}_q$ is a
domain.  Since $d$ clearly is non-zero in $M_q(n)/I^{r+1}_q$ we get the
result. \qed

We now offer an alternative argument in the generic case. We need
this result to prove (\ref{sj-star}) in case $q$ is generic.
Later on, we also obtain (\ref{sj-star}) in the case of a  primitive
$m$th root of unity by ring theoretic methods.

\noindent{\it Proof of Proposition~\ref{5.2} for $q$ generic:} Our
proof relies on   representation 
theory. First of all, for this case it was proved in
\cite{nym} that $M_q(n)$ is a bi-module of a version
$U_q(gl(n,{\mathbb C}))$ of the quantized enveloping algebra of
$gl(n,{\mathbb C})$. Essentially, this version is what results if one
starts from the quantized Serre relations and view the $q$ entering
there as a complex number. Furthermore, it is assumed that $q\neq 0$
and that $q$ is not a root of unity.

The results obtained by \cite{nym} reveal that the same general
picture holds as in the well-known case for $q=1$ \cite{dep}.
Specifically, each $I^{s}_q$ is a $U_q(gl(n,{\mathbb C}))$ sub-bi-module.
Moreover, there is a decomposition
\begin{equation}\label{de-co}
M_{q}(n)=\oplus _{\lambda }W(\lambda )
\end{equation}
as a bi-module. Here,   each $
W(\lambda )$ is an irreducible $U_{q}(gl(n))\times U_{q}(gl(n))$  module. The
highest weight vector in $W_\lambda$ is given by
\begin{equation}\label{v-l}
  w_\lambda=d_1^{a_1}\cdot   d_2^{a_2}\cdot \cdots\cdot d_s^{a_s}
\end{equation}
for $1\leq s\leq n$ and $a_1,a_2,\dots, a_s\in{\mathbb N}\cup\{0\}$.
For each $i\in\{1,\dots,n\}$ let
$\lambda_i=(\underbrace{1,\dots,1}_i,0,\dots,0)$. Then the weight
$\lambda$ of the $w_\lambda$ in (\ref{v-l}) is given by
$\lambda=a_1\lambda_{1}+\dots+a_s\lambda_{s}$. There are no multiplicities.

Furthermore, each $I^s_q$ is invariant. If $W_{q,s}$ denotes the
direct sum of the highest weight modules whose highest weight vectors
are of the form $d_1^{a_1}\cdot d_2^{a_2}\cdot \cdots\cdot
d_s^{a_{\hat s}}$ with ${\hat s}\leq s$, then this is precisely equal
  to $I^s_q\setminus I^{s+1}_q$.
  
  We are now ready to prove that $d$ is regular: $\text{\em Suppose
    that } d\cdot u\in I^{r+1}_q\text{\em . Then } u\in I^{r+1}_q$.
  Assume that $u_q\notin I^{r+1}_q$.  With no loss of generality we
  may assume that $u_q\in W_{q,r}$ and $u_q\neq0$.  Observe that
  $d=d_q$ is a primitive vector for all $q\neq0$.  The 2-sided action
  of the Borel sualgebra ${\mathcal U}^+_{q}(gl(n))$ at $q$ generic
  then preserves the general form of (\ref{de-co}) and hence we may
  assume that $u_q$ is a sum of highest weight vectors of different
  highest weights, and using weight considerations, we may assume that
  $u_q$ is a single highest weight vector. This is then of the form
  (\ref{v-l}) with $s=r$. But we must  still have that $d\cdot
  u_q\in I^{r+1}_q$ and this is a contradiction. \qed
 
\medskip

Below we show that the powers of $d$ can be inverted in a manageable
manner.

\medskip

Let $\pi$ be the natural homomorphism from $A_{n,r}$ to
$M_q^{r+1}(n)$.

\begin{Prop}\label{orepro}
Let $S=\{q^{-i}d^j\mid i,j=0,1,2,\dots\}$. Then $\pi(S)$ is an Ore set of
regular elements in $M^{r+1}_q(n)$.
\end{Prop}

\pf Since $\pi(d)$ is regular, it suffices to prove that $S$ is an Ore
set in $M_q(n)$. Clearly, $S$ is an Ore set in ${\mathbb
  C}\{Z_{i,j}\mid 1\leq i,j\leq r\}$ since $d$ is central in that
algebra. The remaining indeterminates are now added in a suitable
order, i.e. in such a way that $M_q(n)$ is an iterated Ore extension
of ${\mathbb C}\{Z_{i,j}\mid 1\leq i,j\leq r\}$. If $\alpha_{i,j}$
denotes the automorphism corresponding to $Z_{i,j}$, then either
$\alpha_{i,j}(d)=d$ or $\alpha_{i,j}(d)=q^{-1}d$. The result then
follows by \cite[Lemma~1.4]{good}. \qed

In case $q$ is a primitive $m$th root of unity $d^m$ is central by
\cite[Lemma~7.23]{pw} and the localization is then a central localization.

It follows now that $\pi$ induces a homomorphism
$\pi[S^{-1}]:A_{n,r}[S^{-1}]\mapsto M_q^{r+1}[S^{-1}]$, which is onto since
for each $r<i\le n$ or $r<j\le n$ there exists a $t_{i,j}\in A_{n,r}$
such that $d\cdot Z_{i,j}+t_{i,j}$ is a quantum $(r+1)\times(r+1)$
minor.

We can now prove

\begin{Prop}\label{in-pro}The natural homomorphism $\pi:
  A_{n,r}\mapsto M^{r+1}_q(n)$ is injective: 
\end{Prop}

\pf We first give the details for the ``$q$ generic'' case. 

Under the bi-module action of $U_q(gl(n,{\mathbb C}))$ on $M_q(n)$,
the algebra $A_{n,r}$ is invariant under a Borel subalgebra from one
side and under the opposite Borel subalgebra from the other side.
Since $I^{r+1}_q$ is invariant we may argue exactly as in the proof of
Proposition~\ref{5.2}. Thus, if there is a non-zero element $p\in
A_{n,r}$ such that $\pi(p)\in I^{r+1}_q$, there is also a non-zero
highest weight vector $p_h\in A_{n,r}$ such that $\pi(p_h)\in
I^{r+1}_q$. But then $p_h$ is of the form (\ref{v-l}) with $s\geq r+1$
and by looking at leading terms in these quantized determinants, this is
easily seen to be impossible. 

Next suppose $q$ is a primitive $m$th root of unity.

Here all the algebras
\begin{equation*}
A_{n,r}\,,\;A_{n,r}[d^{-1}]\,,\;M_q^{r+1}(n) \mbox{ and }
M_q^{r+1}(n)(d^{-1})
\end{equation*}
are prime affine P.I. algebras.

We recall that for such algebras Krull-dimension, transcendence
degree and G-K-dimension coincide \cite[Proposition 10.6]{mr}. In the
sequel, $\dim$ denotes one of these.

By the definition of transcendence degree any affine
algebra between the original algebra, $B$, and the simple
quotient algebra will have the same dimension as $B$.

By \cite[Proposition 6.5.4]{mr}
\begin{equation*}
\dim A_{n,r}[d^{-1}]\le n^2-(n-r)^2\,,
\end{equation*}
because it is a localization of a $n^2-(n-r)^2$ iterated Ore
extension.

By a recent result \cite{l-r}
\begin{equation*}
\dim M_q^{r+1}(n)\ge n^2-(n-r)^2,
\end{equation*}
and therefore
\begin{equation*}
\dim M_q^{r+1}(n)[d^{-1}]\ge n^2-(n-r)^2\,.
\end{equation*}
If $\Kernel\pi[S^{-1}]\ne 0$ then the Krull-dimension of
$\pi[S^{-1}](A_{n,r}[S^{-1}])$ must be strictly less than
that of $A_{n,r}[S^{-1}]$, therefore $\pi[S^{-1}]$ is injective
and so is $\pi$. \qed

\medskip

\begin{Cor}\label{co1}For $q$ either generic or a primitive root of
  unity, 
\[A_{n,r}[d^{-1}]\simeq M^{r+1}_q(n)[d^{-1}].\]
\end{Cor}

\begin{Cor}\label{com-cor} If $q$ is a primitive root of unity,
$$\deg M_q^{r+1}=\deg A_{n,r}.$$
\end{Cor}

\medskip

\section{Quantized factor algebras of $M_q(n)$.}\label{qfa}    
 
In this section we determine the degree of $A_{n,r}$.

\begin{Thm}\label{com-thm} If $q$ is an odd primitive $m$th
root of unity then $\deg A_{n,r}=m^{nr-r(r+1)/2}$.
\end{Thm}

\pf We fix $r\geq1$ and use induction on $n\geq r$.

If $n=r$  the formula holds by Theorem~\ref{detdeg} in Section~\ref{detsec}. A closer look at the beginning of the proof of that result yields the validity of the formula for $n=r+1$ also.

We view $A_{n,r}$ as an iterated Ore extension:

 \begin{equation}A_{n,r}=A_{n-1,r}[Z_{1,n};\alpha_{1,n}
]\cdots[Z_{r,n};\alpha_{r,n},\delta_{r,n}][Z_{n,1};\alpha_{n,1}]\cdots[Z_{n,r};
\alpha_{n,r},\delta_{n,r}]. 
\end{equation}

The general strategy of the proof is similar to the proof of
Theorem~\ref{detdeg}. We begin by adjoining the variables
$Z_{1,n}\cdots Z_{r,n}$ to $A_{n-1,r}$. Let $A^{(i)}_{n-1,r}$ denote
the algebra obtained by adjoining $Z_{1,n},\dots Z_{i,n}$ so that
$A^{(0)}_{n-1,r}=A_{n-1,r}$. We show that there exist suitable central
elements $c_1,\dots,c_r$ in $A_{n-1,r}$ which behave nicely under the
automorphisms induced by each
of the variables $Z_{1,n}\cdots Z_{r,n}$. This makes it possible to
construct a central element for each $A^{(i)}_{n-1,r}$ which has an
$m$th order orbit under $\alpha_{i+1,n}$. Thus, an application of 
Proposition~\ref{coa} is possible with the conclusion that the degree of
$A^{(r)}_{n-1,r}=A_{n-1,r}[Z_{1,n};\alpha_{1,n}]\cdots[Z_{r,n};\alpha_{r,n},\delta_{r,n}]$
is $m^r$ times the degree of $A_{n-1,r}$.

After that we construct $r$ central elements of $A_{n,r}$. In the general
situation they will have the same shape as those for $A_{n-1,r}$. For
each $Z_{n,i}$ there will in fact be a central element of $A_{n,r}$
which is an inhomogeneous polynomial of degree 1 in that
variable. Thus, Lemma~\ref{cid} implies that the degree does not go up
by adjoining $Z_{n,1},\dots, Z_{n,r}$ to
$A^{(r)}_{n-1,r}$. Having established this, the theorem is proved.

\begin{figure}
\epsffile{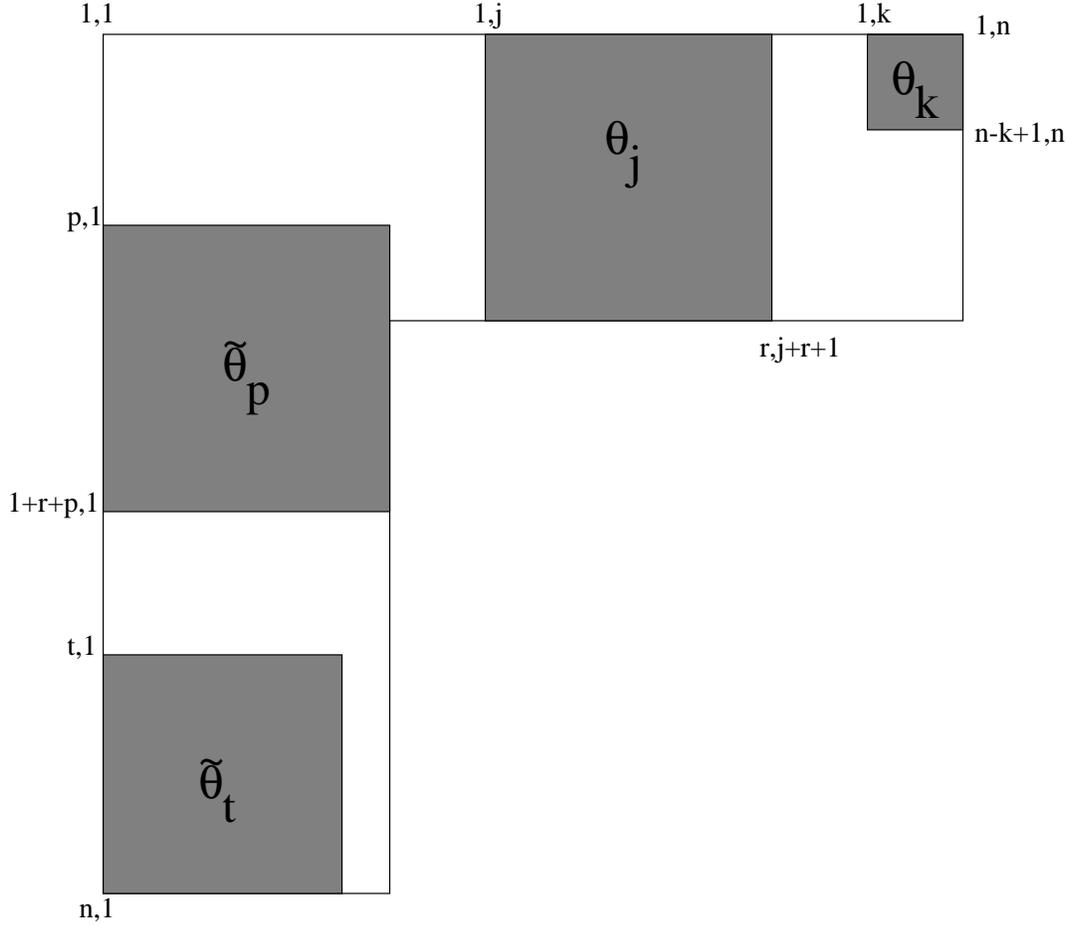}
\caption{\label{f1}The various building blocks.}
\end{figure}

To make the argument clearer we consider the algebra in a diagrammatic
fashion, see Figure~\ref{f1}.

\medskip

We will first consider the cases where $r+1\leq n<3r$. Extend the
previous definitions of $\theta_k, \widetilde \theta_t$ as follows:
For $k> n+1-r$, let $\theta_k$ denote the quantum determinant of the
sub-algebra generated by
 \begin{equation}Z_{1,k},\dots, Z_{1,n},\dots,Z_{n-k+1,k},\dots Z_{n-k+1,n}.
 \end{equation}
 
For $j\leq n+1-r$, $\theta_j$ denotes the quantum determinant (quantum
$r$-minor) of the algebra generated by $Z_{1,j},\dots,
Z_{1,j+r-1},\dots,Z_{r,j},\dots,Z_{r,j+r-1}$.

Elements $\widetilde\theta_j$ are defined analogously for $1\leq j\leq
n$ by interchanging all $Z_{i,j}$ with $Z_{ji}$ in $\theta_j$. It
follows easily that the $\theta_i$ and $\widetilde\theta_i$ are
covariant for all $i=1,\dots,n$.  Notice that
$Z_{i,j}\theta_k=\theta_kZ_{i,j}$ in case $Z_{i,j}$ belongs to the
subalgebra involved in defining $\theta_k$. Also, if $Z_{i,j}$ is not
in the subalgebra used to define $\theta_k$ then
$Z_{i,j}\theta_k=q\theta_kZ_{i,j}$ if $j<k$ and $i\leq
\min\{r,n-k+1\}$. Similarly, $Z_{i,j}\theta_k=q^{-1}\theta_kZ_{i,j}$
in case there exists a $Z_{a,b}$ used in the definition of $\theta_k$
such that $(i,j)=(a,b+x)$ or $(i,j)=(a+x,b)$ for some positive integer
$x$.

In order to make the arguments and ideas as clear as possible, we will
first consider the case where we, starting from ${\mathcal
A}_{r+1,r}$, determine the degree of ${\mathcal A}_{r+2,r}$. We begin
by adjoining the elements $Z_{1,r+2},\dots,Z_{r,r+2}$ to ${\mathcal
A}_{r+1,1}$ in order of increasing first index. Thus,  we
get a sequence of skewpolynomial algebras ${\mathcal A}^{(i)}_{r+1,r}$
for $i=1,\dots,r$. For convenience, let ${\mathcal A}^{(0)}_{r+1,r}=
{\mathcal A}_{r+1,r}$. Now, since ${\mathcal A}_{r+1,r}\subset
M_q(r+1)$ we can use those central elements of the latter that do not
involve $Z_{r+1,r+1}$. In this way we get $r$ central elements:
\begin{equation}c_1=\theta_{r+1}(\widetilde\theta_2^{-1}),\dots, 
c_{r}=\theta_2(\widetilde\theta_{r+1})^{-1},
\end{equation}
where we write $(\widetilde\theta_j)^{-1}$ for $(\widetilde\theta_j)^{m-1}$.

We let $\alpha_j$ denote the automorphism connected with adjoining
$Z_{j,r+2}$ to ${\mathcal A}^{(j-1)}_{r+1,r}$ for $j=1,\dots,r$. Clearly, these automorphisms have order $m$ when acting on the relevant full algebras.

We get 
\begin{equation}\begin{array}{llll}
\alpha_1(c_1)=q^{-1}c_1&\alpha_1(c_2)=q^{-1}c_2&\;\;\qquad\cdots& \alpha_1(c_{r})=q^{-1}c_r\\
&\alpha_2(c_2)=q^{-1}c_2&\;\;\qquad\cdots&\alpha_2(c_{r})=q^{-1}c_{r}\\
\qquad&&\qquad\quad\;\vdots&\qquad\quad\,\,\,\;\;\vdots\\ &&
&\alpha_r(c_{r})=q^{-1}c_{r}.\end{array} 
\end{equation}

Since $c_1$ is central in ${\mathcal A}_{r+1,r}$, and since the length
of the orbit of $\alpha_1$'s action on $c_1$ is $m$, we get that
adjoining $Z_{1,r+2}$ raises the degree by a factor $m$. Next,
$(c_1c_2^{-1})$ is clearly central in ${\mathcal A}^{(1)}_{r+1,r}$ and
$\alpha_2(c_1c_2^{-1})=q^2c_1c_2^{-1}$. Thus, since $m$ is odd, we get
that when we adjoin $Z_{2,r+2}$ to the previously constructed algebra
the degree again goes up by a factor $m$.

Replacing $c_1c_2^{-1}$ by $c_jc_{j+1}^{-1}$ and repeating the
argument, we get  \begin{equation}\deg{\mathcal
A}^{(r)}_{r+1,r}=m^r\deg {\mathcal
A}_{r+1,r}.
\end{equation}

We will now construct $r$ central elements of ${\mathcal
A}_{r+2,r}$. The actual form of these will, in connection with
Lemma~\ref{cid}, imply that adjoining $Z_{r+2,1},\dots,Z_{r+2,r}$ to
${\mathcal A}^{(r)}_{r+1,r}$ does not increase the degree.

First of all, there are $r-1$ central elements coming from $M_q(r+2)$:
\begin{equation}
c_{j}=(\widetilde\theta_{r+3-j})\theta_{j+1}^{-1}\textrm{ for } j=2,\dots,r.
\end{equation}
The remaining central element $c_1$ has got to involve
$Z_{r+2,1}=\widetilde\theta_{r+2}$. We claim that the element
$c_1=\theta_{r+2}\theta_2^{-1}\widetilde\theta_{r+2}(\widetilde\theta_2)^{-1}
\theta_1$ fulfills the requirements. First of all, clearly
$Z_{i,j}c_1=q^{\alpha_{i,j}}c_1Z_{i,j}$ for all $i,j$. To prove
commutativity in all details would involve checking that the five
factors of $c_1$ are in such a balance with each other that the $q^{\alpha_{i,j}}$,
while being the product of five terms of the form $q^{*_k}$ and while
each $q^{*_k}$ depends on the actual form of $(i,j)$, in the end
always equals $q^0$. We leave the somewhat tedious (and somewhat amusing) details of this, as well as similar claims later on, to the reader.

We can now start adjoining the elements $Z_{r+2,i}$. Since for each
$i$, $Z_{r+2,i}$ occurs in the summands of $c_i$ to either the power 1
or 0, Lemma~\ref{cid} applies and the degree remains unchanged.

\medskip

Now suppose $r+1\leq n\leq 2r$.

We have $2r-n+1$ central elements from $M_q(n)$,
\begin{equation}c_{j+1}=\theta_{n-j}(\widetilde\theta_{j+2})^{-1};\qquad
j=n-r-1,\dots,r-1 .\end{equation}

The remaining $n-r-1$ central elements can be chosen as
 \begin{equation}c_{n+1-j}=\theta_1\theta_j(\widetilde\theta_{n-j+2})^{-1}
(\theta_{j-r})^{-1}(\widetilde\theta_{n-j+r+2}),\mbox{ where
}r+1<j\leq n.
 \end{equation}
 
As far as the proof goes, these central elements have the same
properties as the previously constructed. Hence, our strategy applies
and we get that $\deg{\mathcal A}_{n+1,r}=m^r\deg{\mathcal A}_{n,r}$.

\begin{figure}
\epsffile{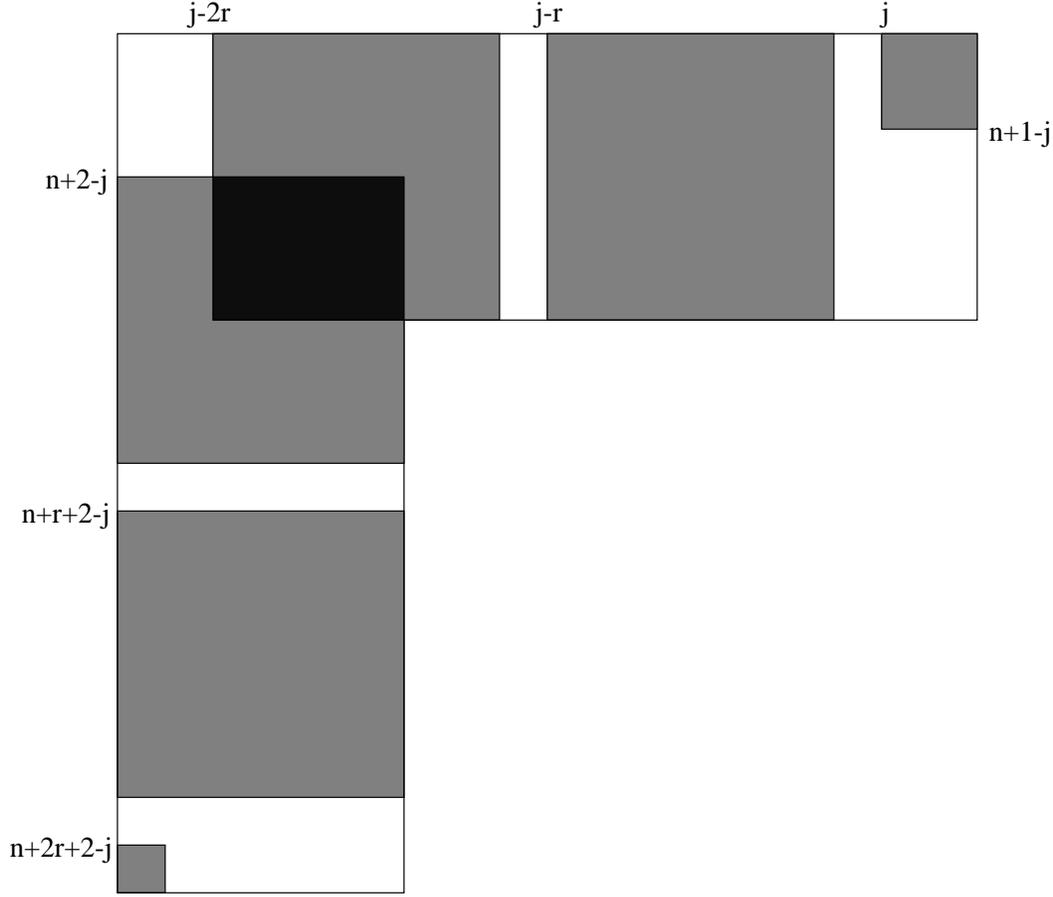}
\caption{Region covered by $\theta_j(\theta_{j-r})^{-1}(\theta_{j-2r})
(\widetilde\theta_{n-j+2})^{-1}(\widetilde\theta_{n-j+2+r})
(\widetilde\theta_{n-j+2+2r})^{-1}$.}
\end{figure}

In the cases $2r< n\leq 3r$ there are no central elements coming from
$M_q(n)$, but there are 2 types of central elements still yielding a
total of $r$ central elements. Specifically for each $j$ with
$n-r+1\leq j\leq 2r+1$ we have the central elements
 \begin{equation}\theta_j(\theta_{j-r})^{-1}\theta_1(\widetilde\theta_{n-j+2})^{-1}
(\widetilde\theta_{n-j+2+r}).
 \end{equation}
The remaining central elements where $j>2r+1$ can be gotten by the
 following:
 \begin{equation}\theta_j(\theta_{j-r})^{-1}(\theta_{j-2r})
(\widetilde\theta_{n-j+2})^{-1}(\widetilde\theta_{n-j+2+r})
(\widetilde\theta_{n-j+2+2r})^{-1}.
\end{equation}
\medskip

\begin{figure}
\epsffile{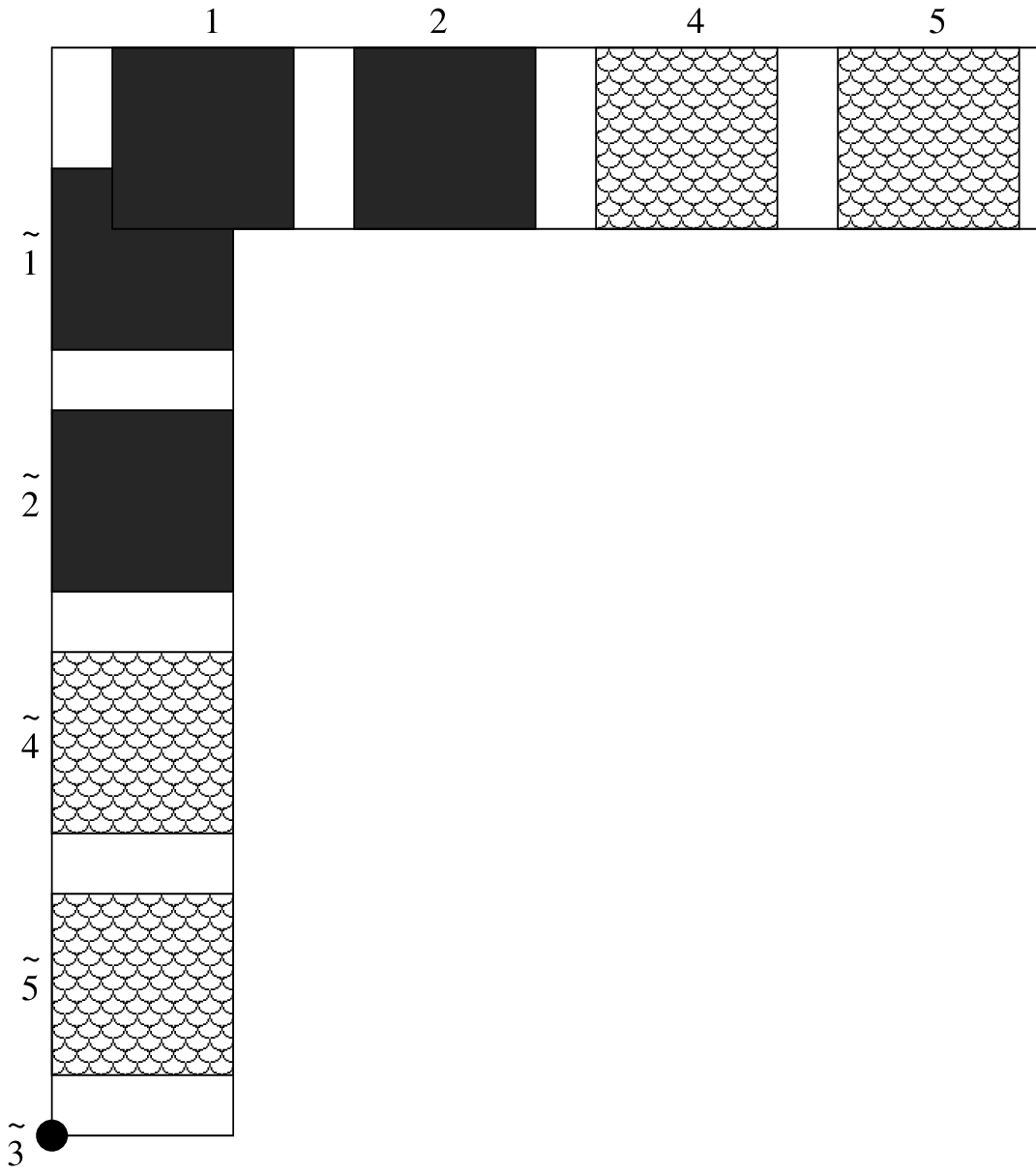}
\caption{\label{fig3} A central element in the general region: $(3)(5)^{-1}(4)(2)^{-1}(1)(\widetilde1)^{-1}(\widetilde2)(\widetilde4)^{-1}(\widetilde5)(\widetilde3)^{-1}$. Compare with Figure 2.}
\end{figure}

Now let us finally comment on the case of a general $n>3r$: Here one
can easily construct $r$ central elements from the previous recipes.
Indeed, observe that each previously constructed central element
contains exactly one pair of factors $\theta,\widetilde\theta$ which
are not full $r\times r$ quantum determinants. If, say, $\theta$ is
$i\times i$ then $\widetilde\theta$ is $(r-i)\times(r-i)$. Even more
precisely, $\theta=\theta_{n-i+1}$ and
$\widetilde\theta=\widetilde\theta_{n-r+i+1}$. Then a number of
factors of the form $(\theta_{n+1-i-k\cdot r})^{\pm1}$ and
$(\widehat\theta_{n-r+i+1-k\cdot r})^{\pm1}$ are inserted for $k\in
{\mathbb N}$ - as long as the resulting indices are positive. Clearly
this procedure continues to work for a general $n$. We omit the finer
details and refer to Figure~\ref{fig3}.  That concludes the proof.
\qed

\medskip

Analogously to  Proposition~\ref{4.6} one gets

\begin{Prop}
  The non-trivial blocks in a block diagonal form of the defining
  matrix $J_A$ of $\overline{A_{n,r}}$ are: $n-1$ matrices of the form
  $\left(\begin{array}{cc}0&1\\-1&0\end{array}\right)$ and \newline
  $nr-\frac{r(r+1)}2 -(n-1)$ matrices of the form
  $\left(\begin{array}{cc}0&2\\-2&0\end{array}\right)$. Let
  $m^\prime=m$ if $m$ is odd and $m^\prime=\frac{m}2$ if $m$ is even.
  Then, in particular, the degree is
  $m^{n-1}(m^\prime)^{nr-\frac{r(r+1)}2 -(n-1)}$.
\end{Prop}

\end{document}